\magnification =1200\hoffset=6mm\voffset=10mm


\normalbaselines
\topskip=18pt
\hsize=120mm \vsize=180mm \parindent=5mm \vskip 15mm

\overfullrule=0mm
\def\det{{\rm det}}
\def\d¡{{\rm d¡}}
\def\tr{{\rm tr}}

\def\adots{\mathinner{\mkern2mu\raise1pt\hbox{.}
\mkern3mu\raise4pt\hbox{.}\mkern1mu\raise7pt\hbox{.}}}

\def\GL{{\rm GL }}
\def\Det{{\rm Det}}

\def\Tr{{\rm Tr}}

\overfullrule=0mm
\def\det{{\rm det}}
\def\d¡{{\rm d¡}}
\def\tr{{\rm tr}}

\def\adots{\mathinner{\mkern2mu\raise1pt\hbox{.}
\mkern3mu\raise4pt\hbox{.}\mkern1mu\raise7pt\hbox{.}}}

\def\GL{{\rm GL }}
\def\Det{{\rm Det}}
\def\Aut{{\rm Aut}}

\def\Spin{{\rm Spin}}

\def\SL{{\rm SL }}

\def\SO{{\rm SO}}
\def\Mat{{\rm Mat}}

\def\codim{{\rm codim}}

\def\Re{{\rm Re}}
\def\Im{{\rm Im}}
\input AMSsym.def
\input AMSsym.tex
\def\tr{{\rm tr}}
\def\Id{{\rm Id}}
\def\Sym{{\rm Sym}}
\def\n{{\rm n }}
\def\Hom{{\rm Hom}}

\def\dim{{\rm dim}}

\input xy.tex
\xyoption{all}

\input xy.tex
\xyoption{all}


\def\boxit#1#2{\setbox1=\hbox{\kern#1{#2}\kern#1}%
\dimen1=\ht1 \advance\dimen1 by #1 \dimen2=\dp1 \advance\dimen2 by #1
\setbox1=\hbox{\vrule height\dimen1 depth\dimen2\box1\vrule}%
\setbox1=\vbox{\hrule\box1\hrule}%
\advance\dimen1 by .4pt \ht1= \dimen1
\advance\dimen2 by .4pt \dp1= \dimen2 \box1\relax}

\bigskip

\bigskip

\centerline{\bf Alg\`ebres de Jordan et th\'eorie des invariants.}

 \bigskip

\centerline{\bf Bruno Blind*}

 \bigskip

\bigskip

\footnote{}{*Institut Elie Cartan, Nancy-Universit\'e, CNRS, INRIA, Boulevard des Aiguillettes, B.P. 239, F- 54506 Vandoeuvre-l\`es-Nancy, France; \hphantom {}blind@iecn.u-nancy.fr}

 \rightline {{\it J'avancerai, je chercherai,}}

 \rightline {{\it jusqu'\`a la fin }}
\bigskip
\rightline {{\it je perdrai courage, je }}

 \rightline {{\it reprendrai de plus loin, je reconna\^\i trai }}

 \rightline {{\it la maison. }}

 \rightline {C. Esteban, La mort \`a distance.}

\bigskip

\noindent{\bf Abstract:} Let $V$ be a simple complex euclidean Jordan algebra of rank three, and denote  by $G$ the subgroup of $\GL(V)$ fixing the determinant of $V$. For $p$ not greater than three, we give a unified description of the invariant  algebras ${\Bbb C}[pV]^G$.
\bigskip

\noindent{\bf R\'esum\'e:} Soit $G$ le sous-groupe de $\GL(V)$ laissant invariant le polyn\^ome d\'eterminant de l'alg\`ebre de Jordan complexe $V$, simple, euclidienne et de rang trois. Dans ce travail nous d\'ecrivons d'une mani\`ere unifi\'ee, pour $p$ inf\'erieur ou \'egal \`a trois, l'alg\`ebre des invariants ${\Bbb C}[pV]^G$.

\bigskip

\bigskip

\bigskip

\noindent{\bf Classification AMS:} 13A50; 15A72; 17C99.

\noindent{\bf Mots clefs:} Th\'eorie classique des invariants, Alg\`ebre de Jordan.

\eject

\beginsection {1. Introduction.}

\bigskip

Si la description de l'alg\`ebre des invariants ${\Bbb C}[pV]^G$ (o\`u $V$ est une repr\'esenta-

\noindent tion de dimension finie d'un groupe alg\'ebrique $G$) est bien connue pour un certain nombre de groupes lin\'eaires classiques ($\SL_n$, $\SO_n,...$), les r\'esultats sont plus rares dans le cas des repr\'esentations de dimension minimale des groupes exceptionnels (voir cependant [Sch.1], [Pop.1] et sa bibliographie ainsi que  les travaux de A.V. Iltyakov cit\'es dans notre bibliographie). 

L'on se propose dans ce travail de donner une description unifi\'ee des g\'en\'era-

\noindent teurs de ${\Bbb C}[pV]^G$ ($p\leq 3$), dans le cas o\`u $V$ est une alg\`ebre de Jordan sur ${\Bbb C}$ simple et euclidienne de rang trois et $G$ le sous-groupe de $\GL(V)$ laissant invariant le d\'eterminant $\det$ de l'alg\`ebre $V$, situation qui comporte en particulier le cas du groupe exceptionnel $G = E_6$ agissant sur l'alg\`ebre d'Albert.

Apr\`es des pr\'eliminaires consacr\'es aux alg\`ebres de Jordan, nous donnons aux paragraphes 3 et 4 des g\'en\'erateurs des alg\`ebres ${\Bbb C}[2V]^G$ et ${\Bbb C}[3V]^G$ en utilisant essentiellement une
m\'ethode g\'eom\'etrique initi\'ee par T. Vust [Vus.]. L'alg\`ebre des invariants ${\Bbb C}[pV]^G$ est, pour $p=2$ une alg\`ebre de polyn\^omes \`a quatre ind\'etermin\'ees; pour $p=3$ et sauf dans le cas $V = \Sym(3,{\Bbb C})$, c'est une alg\`ebre de polyn\^omes \`a onze ind\'etermin\'ees. Ces r\'esultats ne sont pas nouveaux et se trouvent dans la litt\'erature consacr\'ee au sujet (voir [Gor.], [Cia.] et les tables de [Sch.2], de [Ada.-Gol.] et de [Shm.1]), mais nous pensons que notre description unifi\'ee pr\'esente un certain int\'er\^et. Bien entendu, la d\'etermination de ${\Bbb C}[3V]^G$ nous redonne celle de ${\Bbb C}[2V]^G$ (ainsi que celle de ${\Bbb C}[V]^G$).

Dans le cas de l'alg\`ebre de Jordan $V = \Mat(3,{\Bbb C})$, le groupe $G$ n'est pas connexe; nous donnons, au paragraphe 5, une famille g\'en\'eratrice de ${\Bbb C}[3V]^{G^0}$ ainsi que sa s\'erie de Poincar\'e.

Nous terminons, au paragraphe 6, par quelques remarques concernant les alg\`ebres ${\Bbb C}[pV]^G$ pour $p$ sup\'erieur \`a trois.

En appendice, nous avons rappel\'e quelques faits concernant la classification des quatre alg\`ebres de Jordan complexes, simples, euclidiennes et de rang trois.

A certains endroits de ce travail, nous faisons usage du logiciel LIE d\'evelopp\'e par Marc A.A. van Leeuwen, Arjeh M. Cohen et Bert Lisser, nous le d\'esignerons simplement par LIE dans la suite.

Je remercie F. Chargois et O. Hijazi pour leurs critiques et encouragements.

\bigskip

\bigskip

\beginsection {2. Pr\'eliminaires sur les alg\`ebres de Jordan.}

\bigskip

Nous renvoyons \`a [Bra.-Koe.], [Fa.-Ko.] et [Spr.] pour les d\'efinitions et les principaux faits
concernant les alg\`ebres de Jordan.
 
Soit $V$  une des quatre alg\`ebres de Jordan sur ${\Bbb C}$, simples, euclidiennes et de rang $3$: nous  donnons  dans l'appendice une description unifi\'ee de ces quatre alg\`ebres $V_0, V_1, V_2, V_3$ comme l'espace ${\goth {H}}_3({\Bbb F}_{\Bbb C})$ des matrices $(3,3)$ hermitiennes \`a coefficients dans ${\Bbb F}_{\Bbb C}$ (avec ${\Bbb F} = {\Bbb R}, {\Bbb C}, {\Bbb H}, {\Bbb O}$) que l'on munit de la multiplication  de Jordan $x\bullet y = 1/2(xy+yx)$, d'\'el\'ement unit\'e $e = \Id_3$. La forme bilin\'eaire d\'efinie par  $ <x,y> = \tr \ x\bullet y$ est non d\'eg\'en\'er\'ee et  associative:
$$<x\bullet z, y> = <x, y\bullet z>.$$
Nous noterons souvent le complexe $<x,x>$ par $||x||$.

On note par $\det$ le d\'eterminant de $V$, c'est un polyn\^ome irr\'eductible et homog\`ene de degr\'e $3$, sa polarisation nous donne une forme trilin\'eaire $f$ sur $V$ telle que $f(x,x,x) = 6\ \det\ x$. On d\'esigne par  $G$ le sous-groupe  de $\GL(V)$ laissant invariant le d\'eterminant; c'est un sous-groupe du groupe de structure de $V$ (voir [Spr.], Proposition 12.3, page 123) qui contient le groupe $K$ des automorphismes de l'alg\`ebre $V$.

Un \'el\'ement $x$ de $V$ est inversible si $\det\ x \not = 0$, et dans ce cas l'inverse $x^{-1}$ s'\'ecrit sous la forme
$$x^{-1} = {\n(x)\over \det(x)}$$
o\`u $\n$ est une application quadratique de $V$ dans $V$. Cette application est $G$-\'equivariante: 
$$\n(g.x) = (g^{-1})' \n(x)$$
o\`u $g'$ d\'esigne l'adjoint de $g$ par rapport \`a la forme bilin\'eaire $<, >$ (voir [Far.-Kor.] page 148). Par ailleurs, si nous d\'esignons par $\times$ l'application bilin\'eaire sym\'etrique associ\'ee \`a l'application quadratique $\n$: $x\times y = \n(x+y) -\n(x)-\n(y)$, nous avons la relation (voir [Spr.], chapitre 4, formule (5), page 56):

$$f(x, y, z) =<x \times y, z>.$$

L'op\'eration naturelle du groupe $G$ sur $V$ donne lieu \`a une action de $G$ sur l'alg\`ebre ${\Bbb C}[pV]$ des polyn\^omes sur $pV=V\oplus \cdots \oplus V$, et l'objectif de ce travail est de d\'eterminer, dans le cas $p \leq 3$ , l'alg\`ebre des invariants ${\Bbb C}[pV]^G$. On sait (comme cela r\'esulte de la d\'ecomposition de ${\Bbb C}[V]$ en $G$-modules irr\'eductibles, voir [Far.-Kor.]) que ${\Bbb C}[V]^G = {\Bbb C}[\det\ x ]$, l'alg\`ebre \`a une ind\'etermin\'ee
engendr\'ee par le polyn\^ome $\det\ x$,  et \`a partir de la forme trilin\'eaire $f$ et de l'application  $\times $, nous pouvons  facilement construire des invariants du groupe $G$:  par exemple le polyn\^ome $f(\n(x), \n(y),\n(z))$ ainsi que tous ses polaris\'es sont des invariants. Ces invariants ont \'egalement \'et\'e consid\'er\'es par A.V. Iltyakov dans [Ilt.1].

Fixons un syst\`eme complet d'idempotents primitifs orthogonaux $\{e_1, e_2, e_3\}$ de $V$,  on a  la d\'ecomposition de Peirce associ\'ee (voir pour tout ceci les r\'ef\'erences cit\'ees plus haut):

$$V = {\Bbb C} e_1 \bigoplus {\Bbb C} e_2 \bigoplus {\Bbb C} e_2 \bigoplus V_{12} \bigoplus V_{13} \bigoplus V_{23}$$
si nous \'ecrivons la d\'ecomposition d'un \'el\'ement $x$ de $V$ sous la forme:
$$x = x_1e_1+x_2e_2+x_3e_3 + { x}_{12} + {x}_{13} + {x}_{23}$$
nous obtenons facilement les formules suivantes, utiles pour la suite:

$$\n(x)_1 = x_2x_3 - {1\over 2} ||{ x}_{23}||, \ \n(x)_2 = x_1x_3 - {1\over 2} ||{x}_{13}||, \ \n(x)_3 = x_1x_2 - {1\over 2} ||{x}_{12}||$$
et
$$\det \ x = x_1 x_2 x_3 -{1\over 2}( x_1 ||{x}_{23}||+ x_2 ||{ x}_{13}||+ x_3 ||{x}_{12}||)+ 2<{ x}_{12}{ x}_{13},{ x}_{23}>.$$
Si nous \'ecrivons l'\'el\'ement $x$ de $V = {\goth {H}}_3({\Bbb F}_{\Bbb C})$ sous la forme matricielle:
$$ x = \pmatrix{a&{\bar p}&{\bar q}\cr
p&b&r\cr
q&{\bar r}&c}$$
(o\`u $p, q, r$, sont dans ${\Bbb F}_{\Bbb C}$) cette formule devient:
$$\det \ x = a b c - (a\ r {\bar r}+ b\ q {\bar q}+ c\ p {\bar p}) + 2 \Re ((p {\bar q}){\bar r})$$
(o\`u $\Re \ u = 1/2(u+\bar{u})$ pour  un \'el\'ement  $u$ de ${\Bbb F}_{\Bbb C}$). Dans une telle \'ecriture matricielle, on notera souvent les \'el\'ements $x_{12}$, $x_{13}$ et $x_{23}$ par $(p)_{12}$, $(q)_{13}$ et  $(r)_{23}$.

D\'efinissons le groupe $M$, comme le sous-groupe de $K$ fixant les \'el\'ements du syst\`eme d'idempotents  $\{e_1, e_2, e_3\}$. Nous aurons besoin dans la suite de la propri\'et\'e suivante de transitivit\'e du groupe $M$ (nous ne sommes pas parvenus \`a trouver une r\'ef\'erence pour ce r\'esultat, que nous pensons connu). La d\'emonstration proc\'ede cas par cas, pour une preuve unifi\'ee, nous renvoyons \`a la version publi\'ee ([Bli.2]) de ce travail.

\proclaim Proposition 2.1. Soient  $(x_{12},  x_{13}, {x}_{23})$ et $({x'}_{12}, { x'}_{13}, {x'}_{23})$ deux \'el\'ements de $V_{12}\times V_{13}\times V_{23}$ satisfaisant les conditions:
$$\left\{\matrix{ ||{ x}_{ij} || & = & ||{ x'}_{ij} || &\not = 0\cr
<{ x}_{12} \bullet{ x}_{13} , { x}_{23} > & = & <{ x'}_{12} \bullet{ x'}_{13} , {x'}_{23} > & \ \cr}\right.$$ alors il existe un \'el\'ement de $M$ transformant ${ x}_{ij}$ en ${x'}_{ij}$.

L'id\'ee g\'en\'erale est la suivante: par l'action de $M$ nous pouvons nous ramemer au cas o\`u les deux triplets sont de la forme:
 $((a)_{12}, (q)_{13}, (b)_{23})$ et $((a)_{12}, (q')_{13}, (b)_{23})$ (avec $a$ et $b$  deux complexes non nuls); nous devons alors trouver un \'el\'ement de $M$ laissant fixe $(a)_{12}$ et $(b)_{23}$ et transformant $(q)_{13}$ en $(q')_{13}$. La relation 
 $$<(a)_{12}\bullet (q)_{13}, (b)_{23} >  =  <(a)_{12}\bullet (q')_{13}, (b)_{23} >$$ donne alors:
$\Re\ q = \Re\ q'$.
Maintenant, dans les quatre cas, on v\'erifie que le sous-groupe d'isotropie de $(a,b)$ dans $M$ agit transitivement sur la quadrique unit\'e  de $\Im ({\Bbb F}_{\Bbb C}) = \{x\in {\Bbb F}_{\Bbb C}, \Re\ x = 0\}$.

Nous  laissons le d\'etail des cas $V=V_0$ et $V=V_1$ en exercice (noter que pour le cas $V=V_1$, le groupe $M^0$ ne poss\`ede pas cette propri\'et\'e de transitivit\'e), et nous nous concentrons sur les deux autres cas.

Si $V=V_2 =  {\goth {H}}_3({\Bbb H}_{\Bbb C})$, le groupe $M$ est form\'e des matrices diagonales 
$$m_{\lambda,\mu,\nu} = \pmatrix {\lambda&0&0\cr
0&\mu&0\cr
0&0&\nu}$$
(o\`u les trois quaternions complexes sont de norme $1$: $\lambda {\bar \lambda} = \mu {\bar \mu} =\nu {\bar \nu} =1$, voir l'appendice) agissant par:
$$m_{\lambda,\mu,\nu} x = \pmatrix {\lambda&0&0\cr
0&\mu&0\cr
0&0&\nu} x \pmatrix {\bar \lambda&0&0\cr
0&\bar \mu&0\cr
0&0&\bar \nu}.$$
Posons alors
${ x}_{12} =  (p)_{12}, { x}_{13} =  (q)_{13}, { x}_{23} = (r)_{23}$;
l'action d'un \'el\'ement $m_{\lambda,\mu,\nu}$ sur les ${ x}_{ij}$ nous donne:
$$m_{\lambda,\mu,\nu}{ x}_{12} = (\mu \  p\  \bar \lambda)_{12}, \  \  
 m_{\lambda,\mu,\nu}{ x}_{13} = (\nu \ q \ \bar \lambda)_{13}, \  \ m_{\lambda,\mu,\nu}{x}_{23} = (\mu \ r \  \bar \nu)_{23}.$$
en prenant $\mu =1$, $\lambda =  {p \over a}$ et $\nu = { r \over b}$ avec $a$ et $b$ deux complexes tels que $a^2 = p\overline p$ et $ b^2 = r\overline r$, nous pouvons nous ramener au cas o\`u les  deux triplets sont de la forme $((a)_{12}, (q)_{13}, (b)_{23})$ et $((a)_{12}, (q')_{13}, (b)_{23})$ avec $\Re\ q = \Re\ q'$.
Un \'el\'ement  de $M$ laissant fixe $(a)_{12}$ et $(b)_{23}$ est de la forme $m_{\lambda,\lambda,\lambda}$; finalement, nous cherchons un quaternion complexe, $\lambda$, de norme un et tel que
$$ \lambda (\Im\ q ) \overline \lambda = \Im\ q'$$
(o\`u $\Im\ u = 1/2(u-\bar{u})$ pour  un \'el\'ement  $u$ de ${\Bbb F}_{\Bbb C}$). Sachant que $q {\bar q} = q' {\bar q'}$, la  transitivit\'e du groupe $\SO(3, {\Bbb C})$ sur la quadrique $\{(u_1, u_2, u_3) \in {\Bbb C}^3 / u^2_1+u^2_2+u^2_3 = 1\}$ et sur l'ensemble $\{(u_1, u_2, u_3) \in {\Bbb C}^3 / u^2_1+u^2_2+u^2_3 = 0\}\setminus \{(0,0,0)\}$ nous permet de conclure.

Regardons pour finir le cas $V=V_3 =  {\goth {H}}_3({\Bbb O}_{\Bbb C})$; comme rappel\'e dans l'appendice, $M$ est alors le groupe $\Spin(8,{\Bbb C})$: un \'el\'ement $g$ de ce groupe agit par:

$$g.\pmatrix{a&{\bar p}&{\bar q}\cr
p&b&r\cr
q&{\bar r}&c} = \pmatrix {a & \overline{ \rho_+(g)(p)}&  \overline{ \chi (g)(q)} \cr
\rho_+(g)(p) &  b& \rho_-(g)(r) \cr
 \chi (g)(q) &\overline{ \rho_-(g)(r) } & c \cr} $$
o\`u $ \rho_+$, $ \rho_-$ sont les repr\'esentations semi-spinorielles et $\chi $ la repr\'esentation vectorielle de $\Spin(8, {\Bbb C})$. On se ram\`ene \`a nouveau au cas o\`u les deux triplets sont de la forme $((a)_{12}, (q)_{13}, (b)_{23})$ et $((a)_{12}, (q')_{13}, (b)_{23})$ (avec $a$ et $b$ deux complexes non nuls), ce fait nous est garanti ici par les propri\'et\'es de l'action de  $\rho_+ \oplus \rho_-$ sur ${\Bbb C}^8 \oplus {\Bbb C}^8$ (c'est-\`a-dire l'analogue complexe de l'assertion (a) du th\'eor\`eme 14.69 de [Har.], analogue qui se d\'emontre de mani\`ere similaire au cas r\'eel). On a toujours:
$\Re\ q = \Re\ q'$, et le sous groupe d'isotropie de $((a)_{12}, (b)_{23})$ dans $\Spin(8,{\Bbb C})$ s'identifie au groupe complexe $G_2$ agissant  sur $(q)_{13}$ par l'action standard de ce groupe sur ${\Bbb O_{\Bbb C}}$; la transitivit\'e de $G_2$ sur la quadrique $\{(u_1, u_2, \cdots  u_7) \in {\Bbb C}^7 / \sum u_i^2 = 1\}$ de $\Im ({\Bbb O_{\Bbb C}}) = {\Bbb C}^7$ (voir par exemple [Ada.]) et sur l'ensemble $\{(u_1, u_2, \cdots  u_7) \in {\Bbb C}^7 / \sum u_i^2 = 0\}\setminus \{(0,0, \cdots  0)\}$ ach\`event alors la d\'emonstration.

\bigskip
\noindent {\bf Remarque 2.2.} On peut noter que l'application trilin\'eaire:
$$V_{12}\times V_{13}\times V_{23} \mapsto {\Bbb C}$$
$$({ x}_{12}, { x}_{13}, { x}_{23}) \mapsto <{ x}_{12}\bullet { x}_{13} , { x}_{23} >$$
d\'efinit une trialit\'e au sens d'Adams (cf [Ada.] et [Bae.]); le groupe $M$ s'identifie alors au groupe des automorphismes de la trialit\'e. 

\bigskip
Pour terminer ce paragraphe, rappelons un r\'esultat qui nous sera utile pour la suite:

\proclaim Lemme 2.3. L' application de $G$ dans $V$ donn\'ee par:
$g \mapsto g.e$ a pour image $G.e = \{ x\in V /  \  \det x = 1\}$.

Notons $X$ l'ensemble $\{ x\in V /  \  \det x = 1\}$, et soit $x$ un \'el\'ement  de $X$. On sait (proposition VIII 3.5, page 153 de [Far.-Kor.]) que le groupe de structure  de $V$ op\`ere d'une mani\`ere transitive sur l'ensemble des \'el\'ements inversibles de l'alg\`ebre $V$, par cons\'equent il existe un \'el\'ement $g$ de ce groupe tel que 
$$g.e = x.$$
Mais  (voir [Spr.], proposition 1.5, page 11) $\det g.e = \alpha (g) \det e =  \alpha (g)$, donc $g$ est dans $G$; ceci montre que $G$  op\`ere d'une mani\`ere transitive sur $X$, ce qui ach\`eve la d\'emonstration.

\bigskip

\beginsection {3. L'alg\`ebre ${\Bbb C}[2V]^G$.}

\bigskip

Nous d\'eterminons dans ce paragraphe l'alg\`ebre ${\Bbb C}[2V]^G$; pour le cas de l'alg\`ebre des invariants de deux formes quadratiques ternaires (cas $V=V_0$), ce r\'esultat remonte \`a P. Gordan ([Gor.]); les cas $V=V_2$ et $V=V_3$ se trouvent dans les tables de [Sch.2] (voir aussi [Sch.3)]. La d\'emonstration que nous donnons est une illustration de la m\'ethode g\'eom\'etrique de T. Vust ([Vus.]), elle se g\'en\'eralise aux alg\`ebres de Jordan sur ${\Bbb C}$ simples et de rang quelconque.

\proclaim Th\'eor\`eme 3.1. L'alg\`ebre ${\Bbb C}[2V]^G$ est une alg\`ebre de polyn\^omes engendr\'ee par les polarisations de la fonction $\det $:
$$f(x,x,x), \  \   f(y,y,y), \  \  f(x,x,y), \  \  f(y,y,x)$$
avec $(x,y) \in 2V$.

L'id\'ee de la d\'emonstration remonte \`a T. Vust (cf [Vus.]). On identifie 
$2V$ \`a $\Hom ({\Bbb C}^2, V)$: l'\'el\'ement $(x,y)$ de $2V$ est vu comme l'application lin\'eaire: $ (a,b) \mapsto ax+by$, et on consid\`ere l'application $\chi $ de $2V$ dans l'espace $\Sym^3({\Bbb C}^2)$ des polyn\^omes homog\`enes de degr\'e $3$ sur ${\Bbb C}^2$ qui au couple $(x,y)$ associe le polyn\^ome \`a deux variables $a$ et $b$ suivant:
$$(a,b) \mapsto \det (ax+by) = a^3 \det\ x + b^3 \det y +{a^2b\over 2} f(x,x,y) + {ab^2\over 2} f(y,y,x).$$
Comme le morphisme  alg\'ebrique $\chi $ est $G$-invariant, il se factorise suivant le diagramme:

$$\xymatrix{
&2V\ar[ld]_{\pi}\ar[rd]^{\chi}&\\
2V //G\ar@{-->}[rr]_{\Psi}&&{\hbox{Sym}}^3({\Bbb C}^2)}$$
o\`u l'espace $2V //G$ est la vari\'et\'e alg\'ebrique affine irr\'eductible associ\'ee \`a l'alg\`ebre (int\`egre, de type fini) ${\Bbb C}[2V]^G$ et l'application $ \pi: 2V \rightarrow 2V //G$ est le morphisme canonique. 

On va  montrer que $\Psi$ est surjective et birationnelle; il en r\'esultera (voir par exemple [Bri.], lemme 1 page 132) que c'est un isomorphisme alg\'ebrique puisque $\Sym^3({\Bbb C}^2)$ est une vari\'et\'e normale, et le th\'eor\`eme 3.1 s'en suivra.

Montrons que $\chi$, et par cons\'equent $\Psi$, est surjective. On sait qu'un \'el\'ement $P$ (non nul) de $\Sym^3({\Bbb C}^2)$,  peut s'\'ecrire sous la forme:
$$P(a,b) = \lambda b^{3-N}\prod_1^N (a-\alpha_ib)$$
o\`u $\lambda$ est un scalaire non nul, et avec $0\leq N \leq 3$ (un produit vide \'etant \'egal \`a $1$). En posant alors $x = \sum_1^N e_i$ (on prend $x=0$ si $N=0$) et $y = -\sum_1^N\alpha_ie_i + e_{N+1} \cdots +e_3$, il vient 
$$\chi (x,y)  = {1\over \lambda} P$$
et  comme l'image de $\chi$ est ${\Bbb C}^*$ invariante nous obtenons la surjectivit\'e annonc\'ee.

Le morphisme alg\'ebrique $\Psi$ \'etant surjectif, il est dominant; pour montrer qu'il est birationnel il nous suffit (voir par exemple [Vin.], lemme 1 page 252) d'exhiber un sous-ensemble $\cal V$ dense de $\Sym^3({\Bbb C}^2)$ tel que pour tout $v$ de $\cal V$ la fibre $\Psi^{-1}(v)$ soit r\'eduite \`a un point. Pour $(x,y)$ dans $2V$, consid\'erons le polyn\^ome $P(\lambda) = \det (\lambda x-y)$ et son discriminant $\Delta (x,y)$: on constate que pour $\Delta (x,y) \not = 0$ il y a dans la $G$-orbite de $(x,y)$ un \'el\'ement du type $(\lambda e, \sum \lambda_i e_i)$ avec $\lambda$ complexe non nul et les trois coefficients $\lambda_i$ distincts. Posons alors $\cal V = \chi (\cal U)$, $\cal U$ \'etant l'ouvert de Zariski $\{(x,y) \in 2V/ \Delta (x,y) \not = 0 \}$. Soit $(x_0, y_0)$ un \'el\'ement de $\cal U$, on va montrer que si $(x, y) \in 2V$ est tel que 
$$\chi (x_0,y_0) = \chi (x,y) \leqno (*) $$
alors $(x, y)$ et $(x_0, y_0)$ sont dans la m\^eme $G$-orbite, ce qui entra\^\i nera bien que la fibre $\Psi^{-1}(\chi (x_0,y_0))$ est r\'eduite \`a $\pi (x_0,y_0)$. Tout d'abord, on peut supposer que $\det \ x_0 = 1$ sans  restreindre  la g\'en\'eralit\'e du probl\`eme; en se d\'epla\c cant ensuite dans la $G$-orbite de $(x_0, y_0)$, on peut  supposer que $x_0 = e$ et que $y_0$ poss\`ede  trois valeurs propres distinctes. La relation $  (*) $ entra\^\i ne  que $\det \ x = 1$, on peut donc encore, en se d\'epla\c cant cette fois ci dans la $G$-orbite de $(x, y)$ (et quitte \`a changer de notation) remplacer  $(x, y)$ par $(e, y)$; la relation $  (*) $ nous donne en particulier pour tout complexe $a$:
$$\det \ (ae+y_0) =  \det \ (ae+y)$$
et ceci implique que $y$ poss\`ede les m\^emes  valeurs propres distinctes $\lambda_i, i=1, 2, 3$ que $y_0$. On sait alors (voir [Far.-Kor.], proposition VIII.3.2, page 151) que les \'el\'ements $y_0$ et $y$ peuvent s'\'ecrire sous la forme:
$$y_0 =  \lambda_1c_1 +\lambda_2c_2 + \lambda_3c_3, \  \  y =  \lambda_1c'_1 +\lambda_2c'_2 + \lambda_3c'_3$$
o\`u $(c_1, c_2, c_3)$ et  $(c'_1, c'_2, c'_3)$ sont deux rep\`eres de Jordan de $V$, mais l'on sait (voir 
[Bra.-Koe.63], proposition 9.6, page 273 ou [Far.-Kor94], th\'eor\`eme IV.2.5, page 71, la d\'emonstration donn\'ee  pour une alg\`ebre de Jordan r\'eelle s'\'etendant au cas complexe) que le groupe $\Aut(V)$ des automorphismes de $V$ op\`ere transitivement sur les rep\`eres de Jordan. Il r\'esulte de tout ceci que $(x, y)$ et $(x_0, y_0)$ sont bien dans la m\^eme $G$-orbite et ceci ach\`eve la d\'emonstration du th\'eor\`eme.

\bigskip

\bigskip

\noindent {\bf Remarques 3.2.} 

1. Dans le cas de $V=V_1$, on peut noter que les propri\'et\'es de transitivit\'e  des groupes $G$ et $\Aut(V)$ utilis\'ees dans la d\'emonstration sont d\'eja v\'erifi\'ees par leur composante neutre et par cons\'equent  ${\Bbb C}[2V_1]^{G^0} = {\Bbb C}[2V_1]^G$. Nous verrons par la suite que ce r\'esultat ne se g\'en\'eralise pas pour ${\Bbb C}[3V_1]^{G^0}$.

2. Notons \'egalement que pour le cas $V=V_0$, ce r\'esultat nous permet de trouver  la dimension de Krull $d(p)$ de $ {\Bbb C}[pV_0]^G$ pour $p\geq 2$; en effet comme $d(2) = 4$, il r\'esulte de la formule:
$d(2) = \dim \ 2V - \max \dim (G.v)$
que la dimension maximale d'une $G$-orbite $G.v$ dans $2V$ est  $8$, ce qui est la dimension du groupe $G$ dans le cas $V=V_0$, et par cons\'equent $d(p) = 6p-8$. On peut aussi voir directement dans ce cas que le stabilisateur d'un couple de points "g\'en\'eraux" de $V_0$ est fini.

\bigskip

\bigskip

\beginsection {4. Le cas $p = 3$.}

Le but de ce paragraphe est de montrer, en utilisant la m\'ethode g\'eom\'etrique de T. Vust, que les onze polyn\^omes 

$$f(x,x,x), \  \   f(y,y,y),\  \  f(z,z,z),$$
$$f(x,x,y), \  \  f(x,x,z), \  \  f(y,y,x), \  \ 
f(y,y,z),\  \ f(z ,z, x),\  \ f(z,z,y),$$
$$ f(x,y,z), \  \  f(x\times x,y\times y, z\times z).$$
engendrent l'alg\`ebre ${\Bbb C}[3V]^G$; nous noterons d\'esormais ces polyn\^omes $f_1, f_2, \cdots , f_{11}$; on notera que les dix polyn\^omes $f_1, f_2, \cdots , f_{10}$ sont les polarisations de la fonction $\det $.

\noindent On commence par consid\'erer l'analogue du diagramme du paragraphe 2:

$$\xymatrix{
&3V\ar[ld]_{\pi}\ar[rd]^{\chi}&\\
3V //G\ar@{-->}[rr]_{\Psi}&&\Sym^3({\Bbb C}^3)}$$
nous obtenons:

\proclaim Proposition 4.1. Les dix polyn\^omes $f_1, f_2, \cdots , f_{10}$ sont alg\'ebriquement
ind\'epen-
\noindent {\it dants.}

On montre que $\Psi$ est surjective, en prouvant que $\chi$ l'est. Pour ce faire, on remarque
que $\Im \chi = \Im \Psi$ est invariante sous l'action naturelle du groupe $\GL(3,{\Bbb C})$ sur l'espace $\Sym^3({\Bbb C}^3)$, et
donc est une r\'eunion d'orbites de ce groupe; on sait (voir par exemple [Kra.], chapitre 1, paragraphe 7) que  les $\GL(3,{\Bbb C})$ orbites dans $\Sym^3({\Bbb C}^3)$ sont les
orbites de:
$$ a^3, \  \  a^2b, \  \  ab(a+b), \  \  abc,\  \  (a^2-bc)b, \  \  a(a^2-bc)$$
$$b^2c-a^3,\  \  b^2c-a^3-a^2c, \  \  a^3+b^3+c^3+\lambda abc$$
o\`u $\lambda$ est un complexe quelconque. Dans le cas  $V = V_0 = \hbox{Sym}(3, {\Bbb C})$, il n'est pas difficile de v\'erifier que ces repr\'esentants appartiennent \`a $\Im \chi$; le r\'esultat pour les autres cas
en d\'ecoule, au vu des inclusions $V_0 \subset V_1 \subset V_2 \subset V_3$.

Par exemple pour le repr\'esentant $a^3+b^3+c^3+\lambda abc$, on
consid\`ere le triplet $ (e,e_1+\zeta e_2+\zeta^2e_3,z)$ de $3V$ avec $\zeta^3 =1$, et on cherche \`a trouver $z$ tel que 
$$\chi (e,e_1+\zeta
e_2+\zeta^2e_3,z) =  \det (ae+b(e_1+\zeta
e_2+\zeta^2e_3)+cz)=a^3+b^3+c^3+\lambda abc.$$
On a:
$$  \det (ae+b(e_1+\zeta e_2+\zeta^2e_3)+cz) =$$
$$ a^3+b^3+c^3\det z +
ca^2 \tr z + cb^2 \tr ((e_1+\zeta^2 e_2+\zeta e_3) z) + ac^2 \tr \ \n(z) +$$
$$ bc^2 \tr (\n(z) (e_1+\zeta e_2+\zeta^2e_3)) + abc\ \tr ((e_1+\zeta e_2+\zeta^2e_3) \times z).$$
Posons $z =  \pmatrix{
a   & p  &q \cr
p & b &r \cr
q & r & c \cr
}$; les conditions $\tr z = \tr ((e_1+\zeta^2 e_2+\zeta e_3)z) = 0$ et $ \tr ((e_1+\zeta e_2+\zeta^2e_3) \times z) = \lambda $
forment un syst\`eme de Cramer en les inconnues $a, b, c$, et donnent:
$$a =-{\lambda \over 3}, b = -{\lambda \over 3}\zeta^2, c = -{\lambda \over 3}\zeta.$$
Les \'equations
$\tr\ \n(z) = \tr ( \n(z)(e_1+\zeta e_2+\zeta^2e_3)) = 0$ et $\det z =1$ elles, 
sont \'equivalentes \`a:
$$\leqalignno{
&p^2 + q^2 + r^2 = 0 & \cr
& \zeta^2 p^2 + \zeta q^2 + r^2 = 0 & \cr
& \det \ z = 1 & \cr
}$$
ou encore:
$$\leqalignno{
&p^2 + q^2 + r^2 = 0 & \cr
& \zeta^2p^2 + \zeta q^2 + r^2 = 0 & \cr
& {\lambda \over 3}(\zeta p^2 + \zeta^2 q^2 + r^2) + 2pqr =
1+{\lambda^3 \over 27}. & \cr
}$$
Il n'est pas difficile de voir que ce syst\`eme a toujours des solutions: les deux premi\`eres
\'equations permettent par exemple d'exprimer $p^2$ et $q^2$ en fonction de
$r^2$:
$$p^2 = \zeta^2 r^2, \ \ q^2 = \zeta r^2.$$
la derni\`ere condition devenant une \'equation du troisi\`eme degr\'e en l'inconnue $z_{23}$:
$$\lambda r^2 \pm 2 r^3 = 1+{\lambda^3 \over 3^3}$$
qui a toujours des solutions.

\bigskip

\bigskip

Consid\'erons maintenant l'alg\`ebre ${\Bbb C}[3V_0]^G$:  les r\'esultats de Ciamberlini (cf [Cia.]) sur les invariants de trois formes quadratiques ternaires (voir aussi [Shm.1], table 3) entra\^\i nent qu'elle  est engendr\'ee par ses \'el\'ements de degr\'e 3 et 6. Par LIE, les \'el\'ements homog\`enes de degr\'e 3 de 
${\Bbb C}[3V_0]^G$ forment un sous-espace vectoriel de dimension dix, par cons\'equent ce sous-espace est engendr\'e par les polyn\^omes  $f_1, \cdots f_{10}$. De m\^eme, et toujours par LIE, les \'el\'ements homog\`enes de degr\'e 6 forment un sous-espace vectoriel de dimension $56 = 55+1$; il en r\'esulte  que l'alg\`ebre ${\Bbb C}[3V_0]^G$ est engendr\'ee par ${\Bbb C}[f_1, \cdots f_{10}]$ et un \'el\'ement de degr\'e 6; or on voit que le polynome $f_{11}$ n'appartient pas \`a ${\Bbb C}[f_1, \cdots f_{10}]$: dans le cas contraire, en utilisant l'invariance de $f_{11}$ par permutation sur $x,y,z$, on aurait une relation du type:
$$f_{11} = A(f_4f_9 + f_5f_7 + f_6f_8) +B f_{10}^2$$
prenons cette relation pour $x = e_1$ et $y = e_2$, il vient $\n(x) = \n(y) = 0$ et donc $Bf(e_1, e_2, z) = 0$ pour tous les $z$, d'o\`u $B=0$, et par cons\'equent:
$$f_{11} = A(f_4f_9 + f_5f_7 + f_6f_8).$$
On prend \`a nouveau $x = e_1$, et on obtient pour tous les $y$ et $z$: 

$$Af(y,y,e_1)f(e_1,z,z) = 0$$
 d'o\`u $A=0$ ce qui est absurde, puisque l'invariant $f_{11} $ n'est pas identiquement nul. Nous avons donc:

\proclaim Proposition 4.2. Dans le cas $V=V_0$, l'alg\`ebre ${\Bbb C}[3V_0]^G$ est engendr\'ee par les onze  polyn\^omes $f_1, f_2, \cdots , f_{11}$.

\noindent { \bf  Remarque 4.3.} On sait (voir [Shm.1]) que l'alg\`ebre ${\Bbb C}[3V_0]^G$ est d'intersection compl\`ete; d'autre part, dans [Bli.1] nous donnons la s\'erie de Poincar\'e de cette alg\`ebre.

\noindent { \bf  Remarque 4.4.} Il est classique de voir $3V_0$ comme le $\GL(3,{\Bbb C})\times G$ module $({\Bbb C}^3)^* \otimes V_0$ il en r\'esulte une action de $\GL (3, {\Bbb C})$ sur ${\Bbb C}[3V_0]^G$ qui
conserve le degr\'e; LIE nous donne la d\'ecomposition du sous-espace des \'el\'ements homog\`enes de degr\'e 6: on constate la pr\'esence (avec multiplicit\'e un) du $\GL (3, {\Bbb C})$ module $\Sym^2(\wedge^3({\Bbb C}^3))$; il n'est alors pas difficile de v\'erifier que le polyn\^ome:
$$\tilde {f}_{11}= f_{11}  - {2\over 3} (f_4f_9 + f_5f_7 + f_6f_8) + {2\over 3} f_{10}^2$$
est un \'el\'ement non nul de ce module. Le polyn\^ome $f_{11}$ est donc un invariant relatif de  $\GL (3, {\Bbb C})$. Bien entendu cet \'el\'ement, consid\'er\'e dans ${\Bbb C}[3V_i], \  \  0 < i \leq 3$ conserve cette propri\'et\'e d'invariance.

\bigskip

\bigskip

Les r\'esultats pr\'ec\'edents nous conduisent \`a consid\'erer, pour les autres $V_i$, le diagramme commutatif suivant:

$$\xymatrix{
&3V\ar[ld]_{\pi}\ar[rd]^{\tilde{\chi}}&\\
3V //G\ar@{-->}[rr]_{\Psi}&&{\hbox{Sym}}^3({\Bbb C}^3) \bigoplus \Sym^2(\wedge^3({\Bbb C}^3))}$$

\noindent o\`u  l'application ${\tilde{\chi}}$ est donn\'ee par:
$$ {\tilde{\chi}} = (\chi, \  \tilde {f}_{11}).$$

\proclaim Th\'eor\`eme 4.5. Dans les cas  $V= V_1, V_2$ et $V_3$, les onze polyn\^omes $f_1, \cdots f_{11}$ sont alg\'ebriquement ind\'ependants et ils engendrent l'alg\`ebre ${\Bbb C}[3V]^G$.

On commence par montrer que l'image du morphisme $\Psi$ contient l'ouvert form\'e de la r\'eunion des $\GL (3, {\Bbb C})$ orbites des  \'el\'ements de la forme $(a^3+b^3+c^3+\lambda abc, \mu)$, ce qui impliquera en particulier que $\Psi$ est dominant. Par la $\GL (3, {\Bbb C})$ \'equivariance de ${\tilde{\chi}}$, il nous suffit de prouver que les \'el\'ements de la forme $(a^3+b^3+c^3+\lambda abc, \mu)$ appartiennent \`a l'image de $\Psi$, et pour ce faire, on  peut se placer dans le cas o\`u $V = V_1 = \Mat (3,{\Bbb C})$ (puisque $V_1\subset V_2\subset V_3$). On cherche, comme dans la
d\'emonstration de la proposition 4.1,  un triplet $(e, e_1+\zeta e_2+\zeta^2e_3, z)$ (avec $\zeta^3=1$) tel que 
$$ \chi (e, e_1+\zeta e_2+\zeta^2e_3, z) = a^3+b^3+c^3+\lambda abc$$
avec la condition suppl\'ementaire:
$\tilde {f}_{11}(e, e_1+\zeta e_2+\zeta^2e_3, z) =  \mu$.
En posant  $z = (z_{ij})$, on trouve \`a nouveau:
$z_{11} =-{\lambda \over 3}, z_{22} = -{\lambda \over 3}\zeta^2,z_{33} = -{\lambda \over 3}\zeta$.
Les conditions:
$$\tr\ \n(z) = \tr (\n(z) (e_1+\zeta e_2+\zeta^2e_3)) =0$$
$$\tilde {f}_{11}(e, e_1+\zeta e_2+\zeta^2e_3, z) =  \mu$$
sont \'equivalentes \`a un syst\`eme de Cramer en les variables $(z\times z)_{11}$, $ (z\times z)_{22}$ et $(z\times z)_{33} $ :
$$\leqalignno{
&(z\times z)_{11} + (z\times z)_{22} + (z\times z)_{33} = 0, & \cr
& (z\times z)_{11}+ \zeta (z\times z)_{22} + \zeta^2(z\times z)_{33} = 0, & \cr
&(z\times z)_{11} + \zeta^2 (z\times z)_{22} +  \zeta (z\times z)_{33}  =
({2\lambda^2 \over 3}-\mu){1\over 4}. & \cr
}$$
En utilisant les relations:
$$(z\times z)_{11} = 2(z_{22} z_{33} -  z_{23}z_{32}) , (z\times z)_{22} = 2(z_{11} z_{33} -  z_{13}z_{31}), (z\times z)_{33} = 2(z_{22} z_{11} -  z_{21}z_{12})$$
nous obtenons un syst\`eme de Cramer en les variables $z_{21}z_{12}, z_{13}z_{31}, z_{23}z_{32}$ qui donne:
$$z_{12}z_{21} = \zeta^2 {\mu \over 8}, z_{13}z_{31} = \zeta {\mu \over 8},  z_{23}z_{32} = {\mu \over 8}.$$
Compte tenu de ces valeurs, l'\'equation $\det\ z = 1$ devient:
$$z_{12}z_{23}z_{31} + z_{21}z_{13}z_{32} = 1 + {\lambda^3 \over 3} - {\lambda \mu\over 8},$$
mais:
$$(z_{21}z_{13}z_{32}) (z_{12}z_{23}z_{31}) = ({\mu \over 8})^3$$
par cons\'equent  $z_{12}z_{23}z_{31}$ v\'erifie une \'equation du second degr\'e; de tout ceci r\'esulte l'existence de $z$ tel que ${\tilde{\chi}} (e, e_1+\zeta e_2+\zeta^2e_3, z) = (a^3+b^3+c^3+\lambda abc, \mu)$.

\bigskip

Montrons maintenant que le morphisme $\Psi$ est birationnel. En proc\'edant comme dans la d\'emonstration du th\'eor\`eme 3.1, on voit qu'il existe un ouvert de Zariski de $3V$ dont les \'el\'ements $(u, v, w)$ ont dans leur $G$-orbite un triplet de la forme $(\lambda e, \sum\lambda_i e_i, z)$ avec $\lambda$ complexe non nul et les coefficients $\lambda_i $ distincts. Il n'est alors pas difficile de constater que, si deux tels triplets $(\lambda e, \sum\lambda_i e_i, z)$ et $(\mu e, \sum\mu_i e_i, z')$ sont dans la $G$-orbite de $(u, v, w)$ avec $g(\lambda e, \sum\lambda_i e_i, z)=(\mu e, \sum\mu_i e_i, z')$, l'\'el\'ement $g$ est, \`a un facteur multiplicatif pr\`es, dans le groupe $K$ et qu'il permute les \'el\'ements du rep\`ere de Jordan $\{e_1, e_2, e_3\}$; par cons\'equent, si $z =  z_1e_1+  z_2e_2 + z_3e_3 + z_{12}+ z_{13} +  z_{23}$ est l'\'ecriture de $z$ dans la d\'ecomposition de Peirce de $V$, les conditions $<z_{ij}, z_{ij}>\not= 0$ sont ind\'ependantes du choix d'un tel triplet. Soit alors $\cal U$ le sous ensemble dense de $3V$ constitu\'e par les \'el\'ements $(u,v,w)$ dont la $G$-orbite contienne un triplet de la forme $(\lambda e, y, z)$ o\`u $\lambda$ est non nul, $y$ poss\'ede trois valeurs propres distinctes, et les conditions $<z_{ij}, z_{ij}>\not= 0$; on va montrer que pour tout $\xi$ dans ${\tilde \chi} (\cal U)$, la fibre $\Psi^{-1}(\xi)$ est r\'eduite \`a un point. Soient donc deux
triplets $(x,y,z)$ et $(x',y',z')$ de $\cal U$ avec:
$${\tilde \chi} (x,y,z) = {\tilde \chi} (x',y',z')$$
on voit facilement que cette relation est \'equivalente aux deux \'equations:
$ \det (ax+by+cz) = \det (ax'+by'+cz')$ et
$f_{11} (x,y,z) = f_{11}(x',y',z')$. En se d\'epla\c cant  dans l'orbite de $(x,y,z)$ et de
$(x',y',z')$, on peut supposer que $x = x' = e$ et $y = y' = \sum y_ie_i$ (o\`u les 3 valeurs $y_i$
sont distinctes). Si nous explicitons alors les deux \'equations en posant 
$z = z_1e_1+  z_2e_2 + z_3e_3 +  z_{12}+z_{13} + z_{23}$ et
$z' = z'_1e_1+  z'_2e_2 + z'_3e_3 +  z'_{12}+z'_{13} + z'_{23}$,
 nous sommes
conduits \`a examiner des syst\`emes d\'ej\`a rencontr\'es plus haut: 
 l'identit\'e $ \det (ae+b\sum y_ie_i+cz) = \det (ae+b\sum y_ie_i+cz')$ nous donne, en
identifiant les coefficients de $a^2c$, $abc$ et $b^2c$, le syst\`eme de Cramer en les inconnues $z'_i-z_i$:
$$\left\{\matrix{
\hfill (z'_1-z_1) &+&\hfill (z'_2-z_2)&+&\hfill (z'_3-z_3)&= 0, \cr
 (y_2+y_3)(z'_1-z_1) &+& (y_1+y_3)(z'_2-z_2) &+&(y_1+y_2)(z'_3-z_3)&= 0, \cr
 \hfill y_2y_3(z'_1-z_1) &+& \hfill y_1y_3(z'_2-z_2)&+& \hfill y_1y_2(z'_3-z_3)&= 0. \cr
}\right.$$
L'on obtient:
$$z_i = z'_i.$$
De m\^eme, en identifiant les coefficients de $ac^2$,
$bc^2$, et en utilisant  $\tr (\n(y)\times \n(z)) = \tr (\n(y)\times \n(z'))$ (qui est l'\'egalit\'e $f_{11} (e,\sum y_ie_i,z) = f_{11}(e,\sum y_ie_i,z')$),  on trouve un syst\`eme de Cramer en les inconnues $ \n(z')_i-\n(z)_i$ qui nous donne:

$$\n(z)_i = \n(z')_i$$
et donc
$|| z_{ij}||^2 = || z'_{ij}||^2$.
Il vient alors, en mettant ces r\'esultats dans l'\'egalit\'e:
$\det\ z' = \det\  z$:
$$< z_{13}\bullet z_{23},
z_{12}>  = < z'_{13}\bullet z'_{23},
z'_{12}>.$$
Il ne nous reste plus qu'\`a appliquer la proposition 2.1 pour voir qu'il existe un
\'el\'ement du groupe
$M$ qui transforme $z$ en $z'$; les deux
triplets $(x,y,z)$ et $(x',y',z')$ sont donc dans la m\^eme $G$-orbite et par cons\'equent, le morphisme $\Psi $ est bien birationnel.

\bigskip

Il suffit, pour terminer la d\'emonstration du th\'eor\`eme, de montrer que $\Im \Psi$ rencontre toute hypersurface de l'espace ${\hbox{Sym}}^3({\Bbb C}^3) \bigoplus {\Bbb C}$ (voir  le lemme 1 page 132 de [Bri.]). Soit donc  une hypersurface de ${\hbox{Sym}}^3({\Bbb C}^3)
\bigoplus {\Bbb C}$ d'\'equation $H = 0$; on peut supposer que  $H$ contienne effectivement la coordonn\'ee de la composante ${\Bbb C}$, car sinon le r\'esultat est acquis par la d\'emonstration de la proposition 3.1; consid\'erons  alors le coefficient  de $H$ de plus haut degr\'e en cette coordonn\'ee: c'est un polyn\^ome sur ${\hbox{Sym}}^3({\Bbb C}^3)$ qui n'est pas identiquement nul sur l'ouvert form\'e par la r\'eunion des $\GL (3, {\Bbb C})$-orbites d'un
$a^3+b^3+c^3+\lambda abc$, et par cons\'equent il existe bien un \'el\'ement de $\Im \Psi$ dans l'hypersurface consid\'er\'ee.

\bigskip

\noindent { \bf  Remarque 4.6.} La d\'emonstration du th\'eor\`eme 3.3 fait intervenir une propri\'et\'e de transitivit\'e du groupe $M$ qui, dans le cas $V = V_1$, n'est pas v\'erifi\'ee par sa composante neutre. Nous d\'eterminerons, au prochain paragraphe, l'alg\`ebre ${\Bbb C}[3V_1]^{G^0}$.

\bigskip

\noindent { \bf  Remarque 4.7.} Notons \'egalement, pour terminer ce paragraphe, que dans le diagramme suivant:

$$\xymatrix{
&nV\ar[ld]_{\pi}\ar[rd]^{\chi}&\\
nV //G\ar@{-->}[rr]_{\Psi}&& \  \  \  {\overline {\chi(nV)}} \subset \Sym^3({\Bbb C}^n)}$$
(o\`u $n = \dim V$, $\chi (x_1, \cdots x_n)$ d\'esigne le polyn\^ome $\det \sum a_ix_i$ et $\overline {\chi(nV)}$ l'adh\'erence de Zariski  de l'image  de $\chi$ dans $\Sym^3({\Bbb C}^n)$), il est facile de constater que l'application $\Psi$ est birationnelle. En d'autres termes, l'inclusion d'anneaux:
$${\Bbb C}[f(x_i, x_j, x_k)/ \ 1\leq i, j, k \leq n] \rightarrow {\Bbb C}[nV]^G$$
induit un isomorphisme sur les corps des fractions.

\bigskip

\bigskip

\beginsection {5. L'alg\`ebre ${\Bbb C}[3V_1]^{G^0}$.}

Dans tout ce paragraphe, $V$ d\'esignera  l'alg\`ebre de Jordan $\Mat(3, {\Bbb C})$; rappelons que dans ce cas, le groupe $G$ est engendr\'e par sa  composante neutre $G^0$ et la transposition $X \rightarrow {}^tX$; la composante  $G^0$ est le groupe des transformations lin\'eaires:
$(g_{1},g_{2}).X = g_{1}Xg_{2}^{-1}$ o\`u $g_{1}$ et $g_{2}$ sont dans $\SL (3,{\Bbb C})$. 
Consid\'erons  l'\'el\'ement $P$ de ${\Bbb C}[3V]$ suivant:
$$P(x, y, z) = \tr(\n(x) z\ \n(y) x\ \n(z) y).$$

\proclaim Proposition 5.1. Le polyn\^ome $P$ est $G^0$-invariant mais non $G$-invariant.

En effet, la propri\'et\'e de $G^0$-invariance r\'esulte facilement des  formules du type:
$\n((g_{1},g_{2})x)  = g_2 \n(x) g_1^{-1}$. D'autre part, on constate que: $\n(^tx) ={}^t\n(x)$, on en tire:
$$P(^tx, ^ty, ^tz) = \tr({}^t\n(x) ^tz\ {}^t\n(y) ^tx\ {}^t\n(z) ^ty) = \tr ^t(y\ \n(z) x\ \n(y) z\ \n(x)) = $$
$$ \tr (y\  \n(z) x\ \n(y) z\ \n(x)) = \tr  (\n(z) x\ \n(y) z\ \n(x) y) = P(z, y, x).$$
On va montrer que $P(x, y, z) \not = P(z, y, x)$; pour ce faire, on prend $x = e$, $y = a_1e_1 + a_2e_2$ avec $a_1 \not = a_2$ et $a_1a_2 \not = 0$, et on pose $z = (z_{ij})$ et $\n(z) = (Z)_{ij}$. On trouve d'une part:
$P(e, a_1e_1 + a_2e_2, z) = a_1 a_2 (a_1 z_{13} Z_{31} + a_2 z_{23} Z_{32})$
et d'autre part:
$P(z, a_1e_1 + a_2e_2, e) = a_1 a_2 (a_1 Z_{13} z_{31} + a_2 Z_{23} z_{32})$,
comme:
$$Z_{31} = z_{21} z_{32} - z_{22} z_{31}, \   \   \  Z_{13} = z_{12} z_{23} - z_{22} z_{13}$$
$$Z_{32} = z_{12} z_{31} - z_{11} z_{32}, \   \   \  Z_{23} = z_{21} z_{13} - z_{11} z_{23}$$
il vient finalement:
$$P(e, a_1e_1 + a_2e_2, z) - P(z, a_1e_1 + a_2e_2, e) = a_1 a_2(a_1-a_2)(z_{13} z_{32} z_{21} - z_{31} z_{12} z_{23})$$
ce qui est une expression non identiquement nulle.

\proclaim Th\'eor\`eme 5.2. Nous avons:
$${\Bbb C}[3V]^{G^0} = {\Bbb C}[3V]^{G} [P].$$

Les \'el\'ements de ${\Bbb C}[3V]^{G}$ sont ceux de ${\Bbb C}[3V]^{G^0}$ qui sont invariants par la transposition $(x, y, z) \mapsto (^tx, ^ty,^tz)$:
$${\Bbb C}[3V]^{G} = ({\Bbb C}[3V]^{G^0})^{G/G^0}.$$
Comme $G/G^0$ est fini, l'application inclusion:
$${\Bbb C}[3V]^{G} \mapsto {\Bbb C}[3V]^{G^0}$$
fait de ${\Bbb C}[3V]^{G^0}$ un ${\Bbb C}[3V]^{G}$-module de type fini (cf [Kra.] page 111), en particu-

\noindent lier, les polyn\^omes $f_1, f_2, \cdots , f_{11}$ forment un syst\`eme de param\`etres homog\`enes de l'alg\`ebre ${\Bbb C}[3V]^{G^0}$ et la s\'erie de Poincar\'e de cette alg\`ebre s'\'ecrit donc (il est facile de constater que le degr\'e d'un invariant est un multiple de trois):
$$P(t) = {\sum_{i=0}^l a_{3i}t^{3i} \over (1-t^3)^{10}(1-t^6)}.$$
Maintenant, on observe que $ \dim 3V =36 - 3l$ d'apr\`es un r\'esultat de F. Knop (voir [Kno.], on peut aussi consulter  [Pop.2]): nous devons v\'erifier la condition suivante:

$$\codim \{(u,v,w) \in 3V/ \dim G_{(u,v,w)} >0 \} \geq 2$$
o\`u $G_{(u,v,w)}$ est le sous-groupe d'isotropie de $(u,v,w)$. En proc\'edant  comme dans les d\'emonstrations des th\'eor\`eme 3.1 et 4.5, on voit qu'en dehors d'un ensemble alg\'ebrique de codimension sup\'erieur ou \'egal \`a deux, la $G^0$-orbite d'un triplet $(u,v,w)$ contient \`a permutation pr\`es un \'el\'ement de la forme $(\lambda e, \sum \lambda_i e_i, z)$, avec $\lambda$ dans ${\Bbb C}^*$ et o\`u les trois coefficients $\lambda_i$ sont distincts; pour un tel \'el\'ement, le sous-groupe d'isotropie est alors:
$$G^0_{(\lambda e, \sum \lambda_i e_i, z)} = \{ m\in M^0 / mz=z\} = M_z^0.$$
Maintenant, compte tenu de la forme des \'el\'ements du groupe $M^0$, on voit facilement que pour que $M_z^0$ ne soit pas fini, il faut que $z$ soit dans un sous-ensemble alg\'ebrique de dimension cinq. Finalement, le sous-groupe d'isotropie $G^0_{(u,v,w)}$ ne peut-\^etre infini que pour $(u,v,w)$ dans un ensemble de dimension $\dim\ G^0 +1+3+5 = 25$. La condition de F. Knop est donc bien v\'erifi\'ee.

Nous avons donc $3l = 9$; enfin, les anneaux $A(p)$ \'etant de Gorenstein (cf [Pro.], page 562) nous avons $a_0 = a_{9} = 1$ et $a_{9-3i} = a_{3i}$;  finalement:
$$P(t) = {1+a_3t^3 + a_3t^6 + t^9 \over (1-t^3)^{10}(1-t^6)}.$$
Il nous suffit alors de calculer, par utilisation du logiciel LIE, la dimension de l'espace ${\Bbb C}[3V]^{G^0}_3$des \'el\'ements homog\`enes de degr\'e $3$ de ${\Bbb C}[3V]^{G^0}$ pour d\'eterminer le coefficient $a_3$; on trouve: 
$$P(t) = {1+ t^9 \over (1-t^3)^{10}(1-t^6)}.$$
Il existe donc un polyn\^ome $h$ homog\`ene de degr\'e $9$ de ${\Bbb C}[3V]^{G^0}$ tel que ${\Bbb C}[3V]^{G^0}_9 =  {\Bbb C}[3V]^{G}_9\oplus {\Bbb C} h$
et
$${\Bbb C}[3V]^{G^0} = {\Bbb C}[3V]^{G} + {\Bbb C}[3V]^{G} h.$$
Le th\'eor\`eme en r\'esulte facilement.

\bigskip

\noindent { \bf  Remarque 5.3.} Comme me l'a fait remarquer un rapporteur anonyme, on peut \'egalement utiliser les r\'esultats de D.A Shmelkin ([Shm.2], corollaire 3.5.1)  pour \'etablir le th\'eor\`eme 5.2.

\bigskip

\bigskip

\beginsection {6. Quelques remarques sur les alg\`ebres ${\Bbb C}[pV]^{G}$ pour $p>3$.}

\bigskip

On sait (voir [Pro.], th\'eor\`eme 1 de la page 386) que les \'el\'ements de ${\Bbb C}[pV]^{G}$ sont obtenus, pour $p \geq \dim V$, par polarisation de ceux de ${\Bbb C}[(\dim V)V]^{G}$; sauf pour le cas $V=V_0$, nous ne connaissons pas de famille g\'en\'eratrice de cette  alg\`ebre.
Consid\'erons le  cas $V=V_0$: la d\'etermination des invariants (et plus g\'en\'eralement des covariants) de plusieurs formes quadratiques ternaires sous l'action du groupe $G = \SL(3, {\Bbb C})$ est un  probl\`eme qui a \'et\'e abondamment  \'etudi\'e par la th\'eorie classique des invariants; dans ce cas, comme le groupe $G$ est inclus dans $\SL(V)$, le d\'eterminant $\Det (x_1, x_2, \cdots , x_6)$ est un invariant et on sait  (voir [Pro.], th\'eor\`eme 2 de la page 386) que  l'alg\`ebre des invariants de ${\Bbb C}[pV]$ est engendr\'ee par les polaris\'es des invariants de $5$ copies de $V$ et par les d\'eterminants de $6$ copies de $V$. La d\'etermination d'une famille g\'en\'eratrice de ${\Bbb C}[5V]^{G}$ a \'et\'e achev\'ee en 1948 par J.A Todd (pr\'ecisant un travail effectu\'e en 1910 par H.W. Turnbull [Tur]): cet auteur donne dans [Tod.1] et [Tod.2] une telle famille sous forme symbolique. Dans [Bli.1], nous sommes parvenus \`a expliciter ces \'el\'ements \`a partir de la forme trilin\'eaire $f$ et de l'application $\times$, plus pr\'ecisemment, on v\'erifie, en comparant avec la famille de J.A.Todd, que l'alg\`ebre ${\Bbb C}[5V]^{G}$ est engendr\'ee par les polyn\^omes $f(x_i, x_j, x_k)$,  $ f(x_i \times x_j, x_k\times x_l, x_m \times x_n)$ et 
$ f((x_i \times x_j)\times (x_k\times x_l), (x_m \times x_n)\times (x_o\times x_p), x_q)$.
En r\'esum\'e, dans le cas $V=V_0$,  l'alg\`ebre ${\Bbb C}[pV]^G$ est engendr\'ee par les invariants de la forme $f(u,v,w)$ o\`u $u,v,w$ sont des mon\^omes (convenables) en $x_1, x_2, \cdots x_p$ construits avec l'op\'eration $\times$, et par le d\'eterminant $\Det (x_{i_1}, x_{i_2}, \cdots , x_{i_6})$

Nous allons voir que la situation est plus complexe pour les trois autres  $V_i, i>0$. Consid\'erons d'abord l'alg\`ebre ${\Bbb C}[pV]^K$: elle contient la sous-alg\`ebre $\Tr (pV)$ engendr\'ee par les \'el\'ements de la forme $\tr\ u$, o\`u $u$ est un mon\^ome de Jordan en  $x_1, \cdots , x_p$; d'autre part, en utilisant  le logiciel LIE (ou en s'armant de patience), on peut constater qu'il existent, pour certaines valeurs de $p$, des \'el\'ements multilin\'eaires p-altern\'es dans ${\Bbb C}[pV]^K$:  dans le cas de $V_2$ par exemple, il y a un  invariant multilin\'eaire  5-altern\'e (unique \`a un scalaire pr\`es), et dans le cas de $V_3$ il existe une forme multilin\'eaire 9-altern\'ee $K$-invariante. De tels polyn\^omes ne peuvent \^etre dans la sous-alg\`ebre $\Tr (5V)$ (ou dans $\Tr (9V)$): en effet, les \'el\'ements de $\Tr (pV)$ sont de la forme:

$$P = \sum \tr\ u_1 \ \tr\ u_2 \cdots \tr\ u_r$$
o\`u $u_1, \cdots u_r$ sont des mon\^omes de Jordan et avec $1\leq r \leq p$; une telle expression ne peut \^etre p-altern\'ee, car en leur appliquant des op\'erateurs "somme altern\'ee":
$${\cal A}_i: f(x_1, \cdots ,x_i) \mapsto \sum_{\sigma \in {\goth S}_i} \epsilon (\sigma)f(x_{\sigma (1)}, \cdots , x_{\sigma (i)})$$
on d\'etruit des termes dans la somme, par exemple si un  terme $ \tr\  u_1 \ \tr\  u_2 \cdots \tr\ u_r$ contient  $\tr\ [ x_1\bullet (x_2\bullet x_3)]$, il sera d\'etruit par un op\'erateur ${\cal A}_2$ portant sur les variables $x_2$ et $x_3$, puisque $\bullet$ est commutatif.

\noindent {\bf Remarques 6.1. } Ces arguments sont tir\'es de [Ilt.2] (ils ne se trouvent pas dans l'article [Ilt.3] publi\'e sous le m\^eme titre). Signalons que dans [Ilt.2] l'auteur se place dans la situation $V=V_3$ et donne explicitement une forme multilin\'eaire 5-altern\'ee sur $V_3$, invariante sous $K$, malheureusement cet invariant est identiquement nul: pour $V=V_3$, il n' existe  pas de forme multilin\'eaire 5-altern\'ee non nulle et $K$-invariante (v\'erification par LIE). On peut n\'eanmoins voir que la formule donn\'ee dans [Ilt.2] d\'efinit bien une forme multilin\'eaire 5-altern\'ee non nulle et $K$-invariante dans la situation $V=V_2$ (ainsi que pour $V=V_0$ et $V=V_1$). On peut \'egalement v\'erifier que la formule
$$\sum_{\sigma \in {\goth S}_5}  \epsilon (\sigma) \tr \ x_{\sigma(1)} \cdots x_{\sigma(5)}$$
donne, pour les cas $V= V_1, V_2$, la m\^eme (\`a un scalaire multiplicatif pr\`es) forme multilin\'eaire 5-altern\'ee non nulle et $K$-invariante.

\noindent Maintenant, nous avons:

\proclaim Proposition 6.2. Soit $p >1$ un entier, alors l'application:
$${\Bbb C}[pV]^G \longrightarrow {\Bbb C}[(p-1)V]^K$$
$$P(x_1, \cdots, x_p) \mapsto P(e, x_2,\cdots, x_p)$$
est surjective.

En effet, d'apr\`es le principe du transfert (voir par exemple [Gro.]), l'application:
$$({\Bbb C}[G/K] \bigotimes{\Bbb C}[(p-1)V])^G \longrightarrow {\Bbb C}[(p-1)V]^K$$
$$\sum f_i\otimes w_i \mapsto \sum f_i({\bar e})w_i$$
est un isomorphisme d'alg\`ebres. D'autre part, le lemme 2.3 montre que  la restriction:
$${\Bbb C}[V] \rightarrow {\Bbb C}[G.e]$$
$$P \mapsto P_{|_{G.e}}$$
est surjective, puisque $G.e$ est un ensemble Zariski ferm\'e de $V$. Il en r\'esulte que:

$${\Bbb C}[V] \bigotimes{\Bbb C}[(p-1)V] \longrightarrow {\Bbb C}[G/K] \bigotimes{\Bbb C}[(p-1)V] $$
$$P\otimes w \mapsto P_{|_{G.e}}\otimes w$$
est \'egalement surjectif; mais il est classique (voir par exemple la d\'emonstration du lemme 11.3 page 61 dans  [Gro.]) que cette application donne lieu \`a une surjection:

$$({\Bbb C}[V] \bigotimes{\Bbb C}[(p-1)V])^G \longrightarrow ({\Bbb C}[G/K] \bigotimes{\Bbb C}[(p-1)V])^G $$
ce qui ach\`eve la d\'emonstration de la proposition.

\bigskip
D'apr\`es cette proposition, si par exemple $V=V_1$ ou $V_2$, il existe un \'el\'ement $P$ de ${\Bbb C}[6V]^G$, tel que $P(e, x_2,   \cdots x_6)$ soit l'invariant multilin\'eaire  5-altern\'e trouv\'e par LIE et \'evoqu\'e plus haut: comme l'op\'eration $\times$ s'exprime en fonction du produit de Jordan $\bullet$ et que $P(e, x_2,   \cdots x_6)$ n'appartient pas \`a $\Tr (5V)$,  le polyn\^ome $P$ ne peut \^etre combinaison lin\'eaire d'invariants de la forme $f(u,v,w)$ o\`u $u,v,w$ sont des mon\^omes  construits avec  $\times$.

\bigskip
\bigskip

\noindent {\bf Remarque 6.4} Nous voyons donc que la forme des \'el\'ements de l'alg\`ebre ${\Bbb C}[pV]^G$ pour $V=V_i, i>0$ est plus complexe que dans le cas $V=V_0$. Dans cette th\'ematique, les r\'esultats les plus prometteurs nous paraissent \^etre ceux de A. Elduque et A.V. Iltyakov dans [Eld.-Ilt.] (et ceux de [Ilt.3] en ce qui concerne les alg\`ebre ${\Bbb C}[pV]^K$) . Par ailleurs,  en suivant les arguments du paragraphe 4 de [Bli.1], on peut  obtenir les r\'esultats suivants, qui pr\'ecisent un peu ceux de [Eld.-Ilt.]: pour $i=1,2$, l'alg\`ebre ${\Bbb C}[pV_i]^G$ est enti\`ere sur la sous-alg\`ebre engendr\'ee par les invariants $f(x_i, x_j, x_k)$ et $f(x_i \times x_j, x_k\times x_l, x_m \times x_n)$; pour le cas de l'alg\`ebre exceptionnelle $V_3$, il faut adjoindre \`a la sous-alg\`ebre les \'el\'ements de la forme $ f((x_i \times x_j)\times (x_k\times x_l), (x_m \times x_n)\times (x_o\times x_p), x_q)$, mais la d\'emonstration de ces faits, sans \^etre difficile, est longue et assez p\'enible, surtout pour le dernier cas.

\bigskip

\bigskip

\beginsection {Appendice.}

On se propose ici de donner quelques pr\'ecisions sur les  quatre alg\`ebres de Jordan sur ${\Bbb C}$, simples, euclidiennes et de rang $3$: r\'ealisations concr\`etes, groupes $G$, $K$ et $M$. Nos r\'ef\'erences sont  [Bra.-Koe.] et [Spr.] ainsi que [Ada.] et [Har.] pour le cas de l'alg\`ebre exceptionnelle.

On note par ${\Bbb O}$ l'alg\`ebre des octaves de Cayley sur ${\Bbb R}$: 
${\Bbb O} = \{\sum_{i=0}^{i=7} a_if_i \} $
avec $a_i$ dans ${\Bbb R} $, o\`u $f_0$ est l'\'el\'ement unit\'e de l'alg\`ebre (que l'on notera par $1$),  
$f^2_i = -1$ pour $i>0$, $f_if_j = - f_jf_i$ pour $i>0, j>0, i \not=j$, $f_3 = f_1f_2, f_5 = f_1f_4, f_6 = f_2f_4$ et $ f_7 = f_3f_4$. On voit ${\Bbb C}$ dans ${\Bbb O}$ comme la sous-alg\`ebre engendr\'ee par les \'el\'ements $1, f_1$, et  ${\Bbb H}$, l'alg\`ebre des quaternions sur ${\Bbb R}$,  comme la sous-alg\`ebre engendr\'ee par les \'el\'ements $1, f_1, f_2$:
$${\Bbb R} \subset {\Bbb C} \subset  {\Bbb H}\subset {\Bbb O}$$ 
${\Bbb F}$ d\'esignera une de ces  alg\`ebres. L'involution de ${\Bbb O}$:
$u = \sum_{i=0}^{i=7} a_if_i   \mapsto \bar u = a_0 - \sum_{i=1}^{i=7} a_if_i  $
induit les involutions habituelles sur les autres alg\`ebres et ces involutions se prolongent par ${\Bbb C}$ lin\'earit\'e aux complexifi\'ees ${\Bbb F}_{\Bbb C} = {\Bbb F}\bigotimes {\Bbb C}$. Finalement sur $\Mat (3,{\Bbb F}_{\Bbb C})$, on a l'involution:
$$x = (x_{ij})\mapsto ^t\bar x = (\overline {x_{ji}}).$$
Les quatre alg\`ebres de Jordan sur ${\Bbb C}$, simples, euclidiennes et de rang $3$ sont alors les matrices carr\'ees d'ordre $3$, hermitiennes \`a coefficients dans  ${\Bbb F}_{\Bbb C}$:
$${\goth {H}}_3({\Bbb F}_{\Bbb C}) = \{ x \in \Mat (3,{\Bbb F}_{\Bbb C}) / \overline {^tx} = x \}$$
munies du produit de Jordan; on les d\'esignera souvent par $V_0 \subset V_1 \subset V_2 \subset V_3$. Dans cette description, les trois matrices diagonales:
$$ e_1 = \pmatrix {1&0&0\cr
0&0&0\cr
0&0&0},  e_2 = \pmatrix {0&0&0\cr
0&1&0\cr
0&0&0}, e_3 = \pmatrix {0&0&0\cr
0&0&0\cr
0&0&1} $$
forment un syst\`eme complet d'idempotents primitifs orthogonaux, et si 
$$ x = \pmatrix{a&{\bar p}&{\bar q}\cr
p&b&r\cr
q&{\bar r}&c}$$ 
est un \'el\'ement de ${\goth {H}}_3({\Bbb F}_{\Bbb C})$, sa d\'ecomposition de Peirce associ\'ee s'\'ecrit:
$$x =  a e_1 + be_2 + ce_3 +  \pmatrix{0&{\bar p}&0\cr
p&0&0\cr
0&0&0} +  \pmatrix{0&0&{\bar q}\cr
0&0&0\cr
q&0&0} + \pmatrix{0&0&0\cr
0&0&r\cr
0&{\bar r}&0}.$$

Donnons maintenant les groupes $G$, $K$ et $M$ pour chacun des quatre cas.

$\bullet$ $V_0 = \Sym(3, {\Bbb C})$ est l'espace des matrices sym\'etriques sur $\Bbb C$
d'ordre $3$.  Le d\'eterminant est le d\'eterminant usuel, $G$ s'identifie au groupe $\SL (3,{\Bbb C})$ agissant sur $V_0$ par l'action $g.X = gX{}^tg$; le groupe $K$ est alors $\SO (3,{\Bbb C})$. $M$ est le sous-groupe form\'e des matrices   de la forme.
$$\pmatrix {\epsilon&0&0\cr
0&\epsilon'&0\cr
0&0&\epsilon \epsilon'}$$
avec $\epsilon^2 = \epsilon'^2 = 1$.

\bigskip

$\bullet$ $V_1=\Mat(3, {\Bbb C})$ est l'espace des matrices complexes
d'ordre $3$ . Ici encore le d\'eterminant est le d\'eterminant usuel. C'est le seul cas o\`u les groupes $G$ et $K$ ne sont pas connexes: ils sont engendr\'es par leur composante neutre et l'\'el\'ement:
$$X \rightarrow {}^tX.$$
La composante  $G^0$
est le groupe des transformations lin\'eaires $(g_{1},g_{2}).X = g_{1}Xg_{2}^{-1}$ avec $g_{1}$ et
$g_{2}$ dans $\SL (3,{\Bbb C})$, $K_o$ \'etant alors la diagonale $(g,g), g\in \SL (3,{\Bbb C}) $. On trouve que $M_o$ est le sous groupe de $K_o$ des matrices de la forme
$$\pmatrix {\lambda&0&0\cr
0&\mu&0\cr
0&0&{1\over \lambda \mu}}$$
avec $\lambda$, $\mu$ deux complexes non nuls.
\bigskip

$\bullet$ $V_2 = {\goth {H}}_3({\Bbb H}_{\Bbb C})$. Le d\'eterminant est le d\'eterminant de Study (on renvoie  a l'excellent article de synth\`ese  [Asl.]  et \`a sa bibliographie pour tout ce qui touche au d\'eterminant de matrices quaternioniques); $G$ peut \^etre vu comme $\SL (3,{\Bbb H_{\Bbb C}})$ agissant par: $g.X = gX{}^t{\bar g}$, $K$ comme $\{ g \in \SL (3,{\Bbb H_{\Bbb C}}) / g{}^t\bar g = {\Id}_3\}$ et  $M$ comme le groupe des matrices diagonales:

$$\pmatrix {\lambda&0&0\cr
0&\mu&0\cr
0&0&\nu}$$
o\`u $\lambda$, $\mu$ et $\nu$ sont des quaternions complexes de norme $1$.

\bigskip

$\bullet$ Pour l'alg\`ebre de Jordan exceptionnelle $V_3$ on sait (voir par exemple [Spr.]) que
le groupe $G$ (resp. le groupe $K$) est le groupe complexe simplement connexe $E_6$ (resp. le groupe $F_4$).  On constate d'autre part (voir par exemple [Har.] ou [Ada.], ces auteurs traitent le cas de l'alg\`ebre de Jordan exceptionnelle r\'eelle $ {\goth {H}}_3({\Bbb O})$, mais les r\'esultats passent ais\'ement au cas complexe) que $M$ s'identifie au  groupe $\Spin(8, {\Bbb C})$: si $g$ est dans $\Spin(8, {\Bbb C})$, il d\'efinira  un \'el\'ement de $M$ par la formule:

$$\pmatrix{a&{\bar p}&{\bar q}\cr
p&b&r\cr
q&{\bar r}&c} \mapsto \pmatrix {a & \overline{ \rho_+(g)(p)}&  \overline{ \chi (g)(q)} \cr
\rho_+(g)(p) &  b& \rho_-(g)(r) \cr
 \chi (g)(q) &\overline{ \rho_-(g)(r) } & c \cr} $$
(o\`u $ \rho_+$, $ \rho_-$ sont les repr\'esentations semi-spinorielles et $\chi $ la repr\'esentation vectorielle du groupe $\Spin(8, {\Bbb C})$), et r\'eciproquement si $m$ est dans $M$, il laisse stable les espaces $V_{ij}$ de la d\'ecomposition de Peirce et par cons\'equent  son action sur un \'el\'ement  $x$ s'\'ecrit:

$$m. \pmatrix{a&{\bar p}&{\bar q}\cr
p&b&r\cr
q&{\bar r}&c} = \pmatrix {a & \overline{ f_{12}(p)}&  \overline{ f_{13}(q)} \cr
f_{12}(p)(p) &  b& f_{23}(r) \cr
f_{13}(q) &\overline{ f_{23}(r) } & c \cr} $$
en utilisant la trialit\'e, il n'est pas difficile de voir que le triplet $(f_{12}, f_{23}, f_{13})$ est de la forme $(\rho_+(g), \rho_-(g), \chi(g))$ pour un $g$ de $\Spin(8, {\Bbb C})$.

\bigskip

\bigskip

\bigskip

\beginsection {Bibliographie.}

\noindent  [Ada.-Gol.] O.M. Adamovich,E.O. Golovina, {\it Simple linear Lie groups having a free algebra of invariants}, Selecta Math. Soviet. (2) (1983/84),183-220.

\noindent  [Ada.] J. F. Adams, {\it Lectures on exceptional Lie Groups,} The  University of Chicago Press, 1996.

\noindent [Asl.] H. Aslaksen, {\it Quaternionic Determinants,} The Mathe. Intel. 18 (1996) n¡3, 57-65.

\noindent  [Bae.] J. C.  Baez, {\it The Octonions,} Bull. Amer. Math. Soc. (N.S.) 39 (2002) n¡2, 145-205. 

\noindent [Bli.1] B. Blind, {\it La th\'eorie des invariants des formes quadratiques ternaires revisit\'ee,} ArXiv: math.RA/0805.4135.

\noindent [Bli.2] B. Blind, {\it Alg\`ebres de Jordan et th\'eorie des invariants,}  Jour. of Lie Theory, 21, (2011) 123-144.

\noindent  [Bra.-Koe.] H.Braun, M.Koecher, {\it Jordan-Algebren,} Springer Verlag, 1966.

\noindent [Bri.]  M. Brion, {\it Invariants et covariants des groupes r\'eductifs,} in M. Brion, G.W. Schwarz, {\it Th\'eorie des invariants et g\'eom\'etrie des vari\'et\'es quotients,} Paris, Hermann, 2000.

\noindent [Cia.]  C. Ciamberlini, {\it  Sul sistema di tre forme ternarie quadratiche,} Giornale di Battaglini 24 (1886) 141-157.

\noindent [Eld.-Ilt.] A. Elduque, A.V. Iltyakov, {\it On Polynomial Invariants of Exceptional Simple Algebraic Groups.} Canad. J. Math. 51 (1999), 506-522.

\noindent [Far.-Kor.] J. Faraut, A. Kor\'anyi, {\it  Analysis on Symmetric Cones,} Clarendon Press, Oxford 1994.

\noindent [Gor.]  P. Gordan, {\it \"Uber B\"uschel von Kegelschnitten,} Math. Ann 19 (1882) 529-552.

\noindent [Gro.] F.D. Grosshans, {\it Algebraic Homogeneous Spaces and Invariant Theory,} Springer-Verlag LNM 1673, 1997.

\noindent [Har.] F.R. Harvey, {\it Spinors and Calibrations,} Academic Press, 1990.

\noindent [Ilt.1] A.V. Iltyakov, {\it On rational Invariants of the Group $E_6$,} Proc. Ame. Math. Soc. 124 (1996) 3637-3640.

\noindent [Ilt.2] A.V. Iltyakov, {\it Laplace Operator and Polynomial Invariants.} Research Report 97-23, Sydney.

\noindent [Ilt.3] A.V. Iltyakov, {\it Laplace Operator and Polynomial Invariants.} Journal of Algebra 207 (1998), 256-271.

\noindent [Kno.] F. Knop, {\it \"Uber die Glattheit von Quotientenabbildungen,} Manuscripta Math. 56 (1986) 419-427.

\noindent [Kra.] H. Kraft, {\it Geometrische Methoden in der Invariantentheorie,} Vieweg 1985.

\noindent [Pop.1] V.L. Popov, {\it An Analogue of M. Artin's Conjecture on Invariants for Nonassociative Algebras,} Amer. Math. Soc. Transl. 169 (1995), 121-142.

\noindent [Pop.2] V.L. Popov, {\it Groups, Generators, syzygies, and Orbits in Invariant Theory,} AMS Translations of Mathematical Monographs vol. 100, 1992.

\noindent [Pro.] C. Procesi, {\it Lie Groups. An Approach through Invariants and Representations,} Springer-Verlag 2007.

\noindent [Sch.1]  G.W. Schwarz, {\it Invariant theory of $G_2$ and $Spin_7$.} Comment.Math. helvetici 63 (1988), 624-663.

\noindent [Sch.2]  G.W. Schwarz, {\it Representations of Simple Lie Groups with Regular Rings of Invariants,} Inventiones mathematicae 49 (1978) 167-191.

\noindent [Sch.3]  G.W. Schwarz, {\it  When polarizations generate,} Transform. Groups, 12 (2007) n¡4 761-767.

\noindent [Shm.1] D.A. Shmel'kin,{\it On Representations of $\SL_n$ with Algebras of Invariants being Complete Intersections,} Journal of Lie Theory, 11 (2001) n¡1 207-229.

\noindent [Shm.2] D.A. Shmel'kin,{\it On non-connected simple linear groups with a free algebra of invariants,} Izvestiya Math.,60 (1996), 811-856.

\noindent [Spr.]  T.A. Springer, {\it Jordan Algebras and Algebraic Groups,} Springer-Verlag, Ergebn. der Math. und ihrer Grenz., vol 75, 1973.

\noindent [Tod.1]    J.A. Todd, {\it The complete irreducible system of four ternary quadratic,} Proc. Lond. Math. Soc. (2) 52 (1951)   271 - 294.

\noindent [Tod.2]    J.A. Todd, {\it Ternary quadratic types,} Philos. Trans. Roy. Soc. London. Ser. A, 241 (1948)   399 - 456.

\noindent [Tur.]   H.W. Turnbull, {\it Ternary quadratic types,} Proc. London Math. Soc 2 (1910)  81-121.

\noindent [Vin.] E.B. Vinberg, {\it On Invariants of a set of matrices,}  Jour. of Lie Theory, 6, (1996) 249-269.

\noindent [Vus.] T. Vust, {\it Sur la th\'eorie des invariants des groupes classiques,} Ann. Inst. Fourier 26 (1976)  1-31.

\end

\bye